\documentclass[12pt]{article}

\usepackage{amsmath,amssymb}
\begin{document}

\pagestyle{myheadings} \markright{CLASS NUMBERS...}

\def \1{{\bf 1}}
\def \a{{{\mathfrak a}}}
\def \ad{{\rm ad}}
\def \al{\alpha}
\def \ar{{\alpha_r}}
\def \A{{\mathbb A}}
\def \Ad{{\rm Ad}}
\def \Aut{{\rm Aut}}
\def \b{{{\mathfrak b}}}
\def \bs{\backslash}
\def \B{{\cal B}}
\def \c{{\mathfrak c}}
\def \cent{{\rm cent}}
\def \C{{\mathbb C}}
\def \CA{{\cal A}}
\def \CB{{\cal B}}
\def \CC{{\cal C}}
\def \CD{{\cal D}}
\def \CE{{\cal E}}
\def \CF{{\cal F}}
\def \CG{{\cal G}}
\def \CH{{\cal H}}
\def \CHC{{\cal HC}}
\def \CL{{\cal L}}
\def \CM{{\cal M}}
\def \CN{{\cal N}}
\def \CP{{\cal P}}
\def \CQ{{\cal Q}}
\def \CO{{\cal O}}
\def \CS{{\cal S}}
\def \CT{{\cal T}}
\def \CV{{\cal V}}
\def \d{{\mathfrak d}}
\def \det{{\rm det}}
\def \diag{{\rm diag}}
\def \dist{{\rm dist}}
\def \eqn{\begin{eqnarray*}}
\def \endeqn{\end{eqnarray*}}
\def \End{{\rm End}}
\def \F{{\mathbb F}}
\def \Fx{{\mathfrak x}}
\def \FX{{\mathfrak X}}
\def \g{{{\mathfrak g}}}
\def \ga{\gamma}
\def \Ga{\Gamma}
\def \Gal{{\rm Gal}}
\def \h{{{\mathfrak h}}}
\def \Hom{{\rm Hom}}
\def \im{{\rm im}}
\def \Im{{\rm Im}}
\def \Ind{{\rm Ind}}
\def \k{{{\mathfrak k}}}
\def \K{{\cal K}}
\def \l{{\mathfrak l}}
\def \la{\lambda}
\def \lap{\triangle}
\def \li{{\rm li}}
\def \La{\Lambda}
\def \m{{{\mathfrak m}}}
\def \mod{{\rm mod}}
\def \n{{{\mathfrak n}}}
\def \name{\bf}
\def \Mat{{\rm Mat}}
\def \N{\mathbb N}
\def \o{{\mathfrak o}}
\def \ord{{\rm ord}}
\def \O{{\cal O}}
\def \p{{{\mathfrak p}}}
\def \ph{\varphi}
\def \prf{\noindent{\bf Proof: }}
\def \Per{{\rm Per}}
\def \q{{\mathfrak q}}
\def \qed{$ $\newline $\frac{}{}$\hfill {\rm Q.E.D.}\vspace{15pt}\pagebreak[0]}
\def \Q{\mathbb Q}
\def \res{{\rm res}}
\def \R{{\mathbb R}}
\def \Re{{\rm Re \hspace{1pt}}}
\def \r{{\mathfrak r}}
\def \ra{\rightarrow}
\def \rank{{\rm rank}}
\def \supp{{\rm supp}}
\def \Spin{{\rm Spin}}
\def \t{{{\mathfrak t}}}
\def \T{{\mathbb T}}
\def \tr{{\hspace{1pt}\rm tr\hspace{2pt}}}
\def \vol{{\rm vol}}
\def \z{\zeta}
\def \Z{\mathbb Z}
\def \={\ =\ }

\newcommand{\frack}[2]{\genfrac{}{}{0pt}{}{#1}{#2}}
\newcommand{\rez}[1]{\frac{1}{#1}}
\newcommand{\der}[1]{\frac{\partial}{\partial #1}}
\newcommand{\norm}[1]{\parallel #1 \parallel}
\renewcommand{\matrix}[4]{\left( \begin{array}{cc}#1 & #2 \\ #3 & #4 \end{array}
            \right)}

\newcounter{lemma}
\newcounter{corollary}
\newcounter{proposition}
\newcounter{theorex}

\newtheorem{conjecture}{\hspace{-20pt}\stepcounter{lemma} \stepcounter{corollary}
    \stepcounter{proposition}\stepcounter{theorex}Conjecture}[section]
\newtheorem{lemma}{\hspace{-20pt}\stepcounter{conjecture}\stepcounter{corollary}
    \stepcounter{proposition}\stepcounter{theorex}Lemma}[section]
\newtheorem{corollary}{\hspace{-20pt}\stepcounter{conjecture}\stepcounter{lemma}
    \stepcounter{proposition}\stepcounter{theorex}Corollary}[section]
\newtheorem{proposition}{\hspace{-20pt}\stepcounter{conjecture}\stepcounter{lemma}
    \stepcounter{corollary}\stepcounter{theorex}Proposition}[section]

\newtheorem{theorex}{\hspace{-40pt}\stepcounter{conjecture} \stepcounter{lemma}
    \stepcounter{corollary} \stepcounter{proposition}Theorem}[section]
\newenvironment{theorem}{\vspace{30pt}\begin{theorex}}{\end{theorex}\vspace{30pt}}

\title{Class numbers of orders in cubic fields}
\author{Anton Deitmar}

\date{}
\maketitle

{\bf Abstract.} In this paper it is shown that the sum of
class numbers of orders in complex cubic fields obeys an
asymptotic law similar to the prime numbers as the bound
on the regulators tends to infinity. Here only orders are
considered which are maximal at two given primes. This
result extends work of P. Sarnak in the real quadratic
case. It seems to be the first asymptotic result on class
numbers for number fields of degree higher than two.

$$ $$

\tableofcontents

\newpage
\begin{center} {\bf Introduction} \end{center}

$ $

Let $\CD$ be the set of all natural numbers $D\equiv 0,1\ {\rm
mod}\ (4)$, $D$ not a square. For $D\in\CD$ let
$$
\CO_D\=\left\{ \left. \frac{x+y\sqrt{D}}{2}\right| x\equiv
yD\ {\rm mod}\ (2)\right\}.
$$
Then $\CO_D$ is an order in the real quadratic field
$\Q(\sqrt{D})$. As $D$ varies, $\CO_D$ runs through the set of all
orders of real quadratic fields. Let $\epsilon_D>1$ be the
fundamental unit of $\CO_D$ and let $h(D)$ be the class number of
$\CO_D$. C.F. Gau\ss\ noted in \cite{gauss} that, as $x>0$ tends
to infinity,
$$
\sum_{^{D\in\CD}_{D\le x}}h(D) \log(\epsilon_D)\= \frac{\pi^2
x^{\frac{3}{2}}}{18\zeta(3)} +O(x\log x).
$$
This was confirmed later by C.L Siegel in \cite{siegel}.

Note that $\log\epsilon_D$ equals the regulator $R(\CO_D)$ of the
order $\CO_D$. For a long time it was believed to be impossible to
separate the class number and the regulator. However, in 1981 P.
Sarnak showed \cite{sarnak}, using the trace formula, that
$$
\sum_{^{D\in \CD}_{\epsilon_D \le x}}h(D) \= \li(x^2) +
O(x^{\frac{3}{2}}(\log x)^2),
$$
where
$$
\li(x)\=\int_2^\infty \rez{\log t} dt
$$
is the integral logarithm as it appears in the prime number
theorem (modulo the Riemann hypothesis).  Let
$r(\CO_D)=e^{2R(\CO_D)}$. Then Sarnak's results can be rewritten
as
$$
\sum_{^{D\in \CD}_{r(\CO_D) \le x}}h(D) \= \li(x) +
O\left(\frac{x^{\frac{3}{4}}(\log x)^2}4\right).
$$
Sarnak established this result by identifying the regulators with
lengths of closed geodesics of the modular curve $H/SL_2(\Z)$
(Theorem 3.1 there) and by using the geodesic prime number
theorem for this Riemannian surface.

In this paper we consider the class numbers of orders in
cubic fields. To our knowledge, the following is the first
asymptotic result on class numbers for number fields of
higher degree.

A cubic field is either totally real, or it has one real
and two complex embeddings. In the second case, we call it
a complex cubic field. Let $p$ and $q$ be two distinct
prime numbers and let $C(p,q)$ denote the set of all
complex cubic number fields $F$, in which neither $p$ nor
$q$ is decomposed; i.e. a complex cubic field $F$ lies in
$C(p,q)$ if and only if there is only one prime of $F$
above $p$ and only one above $q$. Next let $O(p,q)$ be the
set of all orders $\CO$ of fields in $C(p,q)$, which are
maximal at $p$ and at $q$. Let $h(\CO)$ denote the class
number of such an order $\CO$. Let $\la(\CO)$ be $1$, if
$p$ and $q$ are ramified in $F$, let it be $3$ if only one
of them is ramified and $9$ if neither is ramified. For
$x>0$ set
$$
\pi_{p,q}(x)\= \sum_{^{\CO\in O(p,q)}_{r(\CO)\le
x}}h(\CO)\la(\CO),
$$
where $r(\CO)=e^{3R(\CO)}$, and $R(\CO)$ is the regulator of the
order $\CO$. The main result of this paper is
$$
\pi_{p,q}(x)\ \sim\ \frac{x}{\log x}
$$
as $x$ tends to infinity. More precisely we have
$$
\pi_{p,q}(x)\= \li(x) + O\left( \frac{x^{\frac{3}{4}}}{\log
x}\right).
$$
The same conclusion holds if the set $\{ p,q\}$ is replaced by an
arbitrary finite set of primes $S$ with at least two elements. It
is this latter version which we give in the text.

This paper grew out of a conversation with Peter Sarnak
who urged me to look for applications of the theory
developed in \cite{geom}. I thank him warmly for drawing
my attention into this direction. I also thank Nigel Byott
for some useful discussions.

\newpage
\section{The main theorem}
Let $\CO$ be an order in a number field $F$. Let $I(\CO)$
be the set of all finitely generated $\CO$-submodules of
$F$. According to the Jordan-Zassenhaus Theorem
\cite{reiner}, the set of isomorphism classes $[I(\CO)]$
of elements of $I(\CO)$ is finite. Let $h(\CO)$ be the
cardinality of the set $[I(\CO)]$, called the {\it class
number} of $\CO$.

The multiplicative group of invertible elements $F^\times$ of $F$
acts on $I(\CO)$ by multiplication: $\la .M=\la M=M\la$. Since
$M$ and $M\la$ are isomorphic as $\CO$-modules we get a map
$I(\CO)/F^\times\ra [I(\CO)]$ mapping $MF^\times$ to its class
$[M]$. We claim that this map is a bijection. It is clearly
surjective. So let $M$ and $N$ be elements in $I(\CO)$ which are
isomorphic. Fix an isomorphism $T:M\ra N$. Then $T$ extends to an
$F$-isomorphism $T_F : F=FM\ra FN= F$, so there is $\alpha\in
F^\times$ such that $T_F$ is just multiplication by $\alpha$. This
implies the claim. It follows that the class number $h(\CO)$
equals the cardinality of $I(\CO)/F^\times$.

The class number of the maximal order is also called the class
number of the field. In general the class number of $\CO$ will be
larger than that of $F$.

A number field $F$ of degree $3$ over the rationals is also
called a {\it cubic field}. Consider the set of embeddings of a
given cubic field into $\C$. This set has three elements and is
permuted by complex conjugation. It follows that either all three
embeddings are fixed by complex conjugation, i.e. they are real
in which case the field is called {\it totally real} or two of
them are swapped by complex conjugation and the third is fixed,
i.e. one is real and the other two are complex, in which case the
field is called {\it complex}.

Let $F$ be a complex cubic field. Since the order of the
automorphism group of $F$ divides the degree, which is $3$, it is
either $1$ or $3$. In the latter case the extension would be
galois, so the Galois group would act transitively on the set of
embeddings of $F$ into the complex numbers. In particular, all
embeddings would be either real or complex. Since this is not the
case it follows that the automorphism group of a complex cubic
field is trivial.

Let $F$ be a complex cubic field and let $\CO$ be an order of
$F$. Then the group of units satisfies
$$
\CO^\times\= \pm \epsilon^\Z
$$
for some base unit $\epsilon$. The image of the base unit
under the real embedding will be of the form $\pm e^{\pm
R(\CO)}$, where $R(\CO)$ is the {\it regulator} of $\CO$
(see \cite{neukirch}). Set
$$
r(\CO)\= e^{3R(\CO)}.
$$

A prime number $p$ is called {\it non-decomposed} in $F$ if there
is only one place in $F$ lying above $p$. Fix a finite set $S$ of
prime numbers with at least two elements and let $C(S)$ be the
set of all complex cubic fields $F$ such that all primes $p$ in
$S$ are non-decomposed in $F$. For $F\in C(S)$ let $O_F(S)$ be
the set of all orders $\CO$ in $F$ which are maximal at all $p\in
S$, i.e. are such that the completion $\CO_p=\CO\otimes_\Z\Z_p$ is
the maximal order of the field $F_p=F\otimes_\Q\Q_p$ for all
$p\in S$. Let $O(S)$ be the union of all $O_F(S)$ where $F$
ranges over $C(S)$.

Let $S_i(F)$ be the set of all $p\in S$ such that $p$ is
inert in $F$. Define
\eqn
\la_S(F) &=& d^{|S_i(F)|}\\
&=& \prod_{p\in S} f_p(F),
\endeqn
where $f_p(F)$ is the inertia degree of $p$ in $F$. For an
order $\CO$ in $F$ let
$$
\la_S(\CO)\= \la_S(F).
$$

The following is our main theorem.

\begin{theorem}\label{main}
For $x>0$ let
$$
\pi_{S}(x) \= \sum_{^{\CO\in O(S)}_{r(\CO)\le
x}}h(\CO)\la_S(\CO).
$$
Then, as $x\ra\infty$ we have
$$
\pi_{S}(x)\ \sim\ \frac{x}{\log x}.
$$
More sharply,
$$
\pi_{S}(x)\= li(x)+O\left(\frac{x^{\frac{3}{4}}}{\log
x}\right),
$$
as $x\ra\infty$ where $li(x)=\int_2^\infty \rez{\log t}dt$
is the integral logarithm as in the usual prime number
theorem.
\end{theorem}

As an immediate consequence we get the following corollary.

\begin{corollary}
Let
$$
\tilde\pi_S(x)\=\sum_{\frack{\CO\in O(S)}{r(\CO)\le x}}
h(\CO),
$$
then
$$
\limsup_{x\ra\infty}\frac{\tilde\pi_S \log x}x\ \le\ 1,
$$
and
$$
\liminf_{x\ra\infty}\frac{\tilde\pi_S \log x}x\ \ge\ \frac
1{3^{|S|}}.
$$
\end{corollary}

\prf For any $\CO$ we have $1\le \la_S(\CO)\le 3^{|S|}$,
which implies the corollary.
\qed

The theorem will be proved in the following sections.

To clarify the range of the theorem we note:

\begin{proposition}\label{nigel}
For every complex cubic field $F$ there are infinitely
many primes $p$  which are non-decomposed in $F$. Every
order $\CO$ in $F$ is maximal at almost all primes $p$.
\end{proposition}

\prf The second assertion is well known and so it simply remains
to prove the first. Let $F$ be a complex cubic field. We will show
that there are infinitely many primes $p$ which are
non-decomposed in $F$. Let $E/\Q$ be the galois hull of $F$. Then
$\Gal(E/\Q)$ is the permutation group $S_3$ in three letters. Let
$\sigma$ be the generator of $\Gal(E/F)\cong\Z/2\Z$ and let $p$ be
a prime which is unramified in $E$ and such that there is a prime
$\p_1$ of $E$ over $p$ with Frobenius $\tau$ of order $3$. By the
Tchebotarev density theorem there are infinitely many such $p$.
Then the stabilizer of $\p_1$ in $S_3$ is $\langle\tau\rangle$ of
order $3$, so $p\CO_E=\p_1\p_2$, say. As $\sigma$ does not
stabilize $\p_1$, it interchanges $\p_1$ and $\p_2$, so $\p_1\cap
F=\p_2\cap F$ is a prime in $\CO_F$, i.e. $p$ is non-decomposed
in $F$. \qed

\section{Division algebras of prime degree}
Let $d$ be a prime $>2$. Let $M(\Q)$ denote a division algebra of
degree $d$ over $\Q$. The dimension of $M(\Q)$ is $d^2$. Fix a
maximal order $M(\Z)\subset M(\Q)$; for any commutative ring with
unit $R$ we will write $M(R)\= M(\Z)\otimes_\Z R$. Further, let
$M(R)^\times$ denote the multiplicative group of invertible
elements in $M(R)$. We say that $M(\Q)$ {\it splits} over a prime
$p$ if $M(\Q_p)$ is not a division algebra. Since $d$ is a prime
we then get an isomorphism $M(\Q_p)\cong\Mat_d(\Q_p)$ (see
\cite{pierce}). Further, since there are no associative division
algebras over the reals of degree $d$ it follows that
$M(\R)\cong\Mat_d(\R)$. For any ring $R$ the reduced norm induces
a map $\det:M(R)\ra R$. Let $\CG(R)=\{ x\in M(R) | \det(x)=1\}$.
Then $\CG$ is a simple linear algebraic group over $\Z$. Let
$G=\CG(\R)$; then $G$ is isomorphic to $SL_d(\R)$. Let $S$ be the
set of places of $\Q$ where $M(\Q)$ does not split. The set $S$
always is finite, has at least two elements, and contains only
finite places. For any prime $p$ we have that $M(\Z_p)$ is a
maximal $\Z_p$-order of $M(\Q_p)$.

Let $F/\Q$ be a finite field extension which embeds into
$M(\Q)$. By the Skolem-Noether Theorem (\cite{reiner}, Thm
7.21), any two embeddings $\sigma_1,\sigma_2 : F\ra M(\Q)$
are conjugate by $M(\Q)^\times$, i.e. there is a $u\in
M(\Q)^\times$ such that $\sigma_2(x)=u\sigma_1(x)u^{-1}$
for any $x\in F$.

A prime $p$ is called {\it non-decomposed} in the
extension $F/\Q$ if there is only one place $v$ of $F$
lying above $p$.

\begin{lemma}
Let $F/\Q$ be a nontrivial finite field extension. Then
$F$ embeds into $M(\Q)$ if and only if $[F:\Q]=d$ and
every prime $p$ in $S$ is non-decomposed in $F$.
\end{lemma}

\prf Assume that $F$ embeds into $M(\Q)$; then $[F:\Q]$ divides
$d$, which is a prime, so that $[F:\Q]=d$. It follows that $F$ is
a maximal subfield of $M(\Q)$. By Proposition 13.3 in
\cite{pierce} it follows that $F$ is a splitting field of
$M(\Q)$. Thus, by Theorem 32.15 of \cite{reiner} it follows that
for every $p\in S$ and any place $v$ of $F$ over $p$ we have $[F_v
:\Q_p]=d$, which by the degree formula implies that there is only
one $v$ over $p$. Conversely, assume that there is only one prime
above each $p\in S$. Let $p\in S$ and let $e$ be the ramification
index of $p$ in $F$ and $f$ the inertia degree. Then
$d=[F:\Q]=ef=[F_v:Q_p]$ by the global, respectively local, degree
formula. Since this holds for any $p\in S$, \cite{reiner},
Theorem 32.15 implies that $F$ is a splitting field for $M(\Q)$.
By Proposition 13.3 of \cite{pierce} we infer that $F$ embeds
into $M(\Q)$. \qed

Let $F/\Q$ be a field extension of degree $d$ which embeds
into $M(\Q)$. Then for any embedding $\sigma :F\ra M(\Q)$
the set
$$
\CO_\sigma \= \sigma^{-1}(\sigma(F)\cap M(\Z))
$$
is an order of $F$. For $p\in S$ let $v_p$ denote the
unique place of $F$ over $p$.

\begin{lemma}\label{1.2}
Let $\sigma :F\ra M(\Q)$ be an embedding of the field $F$. For
any $p\in S$ the completion $\CO_{\sigma, v_p}$ is a maximal order
in $F_{v_p}$. Conversely, let $\CO\subset F$ be an order such that
for any $p\in S$ the completion $\CO_{v_p}$ is maximal. Then
there is an embedding $\sigma :F\ra M(\Q)$, such that $\CO
=\CO_\sigma$.
\end{lemma}

\prf Let $p\in S$; then $M(\Z_p)$ is a maximal $\Z_p$-order in the
division algebra $M(\Q_p)$, and hence by Theorem 12.8 of
\cite{reiner} it coincides with the integral closure of $\Z_p$ in
$M(\Q_p)$. Therefore
$\CO_{\sigma,v_p}=\sigma^{-1}(\sigma(F_{v_p}\cap M(\Z_p))$ is the
integral closure of $\Z_p$ in $F_{v_p}$, which is the maximal
order of $F_{v_p}$.

For the converse, let $\CO$ be an order of $F$ such that the
completion $\CO_{v_p}$ is maximal for each $p\in S$. Fix an
embedding $F\hookrightarrow M(\Q)$ and consider $F$ as a subfield
of $M(\Q)$. For any $u\in M(\Q)^\times$ let $\CO_u=F\cap
u^{-1}M(\Z)u$. We will show that there is a $u\in M(\Q)^\times$
such that $\CO=\CO_u$. This will prove the proposition since one
can then take $\sigma$ to be the conjugation by $u$.

Let $\CO_1=F\cap M(\Z)$. Since $\CO$ and $\CO_1$ are orders, they
both are maximal at all but finitely many places. So there is a
finite set of primes $T$ with $T\cap S=\emptyset$ and such that
with $T_F$ denoting the set of places of $F$ lying over $T$ we
have that for any place $v$ of $F$ with $v\notin T_F$ the
completion $\CO_v$ is maximal and equals $\CO_{1,v}$. Let $p\in
T$ and fix an isomorphism $M(\Q_p)\ra \Mat_d(\Q_p)$. Fix a
$\Z_p$-basis of $\CO\otimes_\Z \Z_p$. This basis then induces an
embedding $\sigma_p : F\otimes_\Q \Q_p\ra \Mat_d(\Q_p)=M(\Q_p)$,
such that
$$
\sigma_p^{-1} ( \sigma_p(F\otimes\Q_p)\cap\Mat_d(\Z_p))\=
\CO\otimes\Z_p.
$$
By the Noether-Skolem Theorem \cite{reiner} 7.21 there is a
$\tilde{u}_p\in M(\Q_p^\times)$ such that
$\CO\otimes\Z_p=F_p\cap\tilde{u}_p^{-1}M(\Z_p)\tilde{u}_p$. For
$p\notin T$ set $\tilde{u}_p=1$ and let
$\tilde{u}=(\tilde{u}_p)\in M(\A_{fin})$, where $\A_{fin}$ is the
ring of the finite adeles over $\Q$, i.e. the restricted product
over the local fields $\Q_p$, where $p$ ranges over the primes.
By strong approximation $M(\Q)^\times$ is dense in $M(\A_{fin})$
and so there is $u\in M(\Q)^\times$ such that $M(\hat{\Z})
u=M(\hat{\Z})\tilde{u}$, where $\hat{\Z}=\prod_p\Z_p$. It follows
that $\CO=\CO_u$ for this $u$.
\qed

Let $F$ be a field extension of $\Q$ of degree $d$ which
embeds into $M(\Q)$. Let $\CO$ be an order of $F$, which
is maximal at each place in $S$. By Lemma \ref{1.2} we
know that there is an embedding $\sigma$ of $F$ into
$M(\Q)$ such that $\CO =\CO_\sigma$. Let $u\in
M(\Z)^\times$ and let $^u\sigma = u\sigma u^{-1}$. Then
$\CO_{^u\sigma}=\CO_\sigma$, so the group $M(\Z)^\times$
acts on the set $\Sigma(\CO)$ of all $\sigma$ with
$\CO=\CO_\sigma$.

\begin{lemma}
The quotient $\Sigma(\CO)/M(\Z)^\times$ is finite and has
cardinality equal to the product $h(\CO)\la_S(\CO)$.
\end{lemma}

\prf Fix an embedding $F\ra M(\Q)$ and consider $F$ as a
subfield of $M(\Q)$ such that $\CO=F\cap M(\Z)$. For $u\in
M(\Q)^\times$ let
$$
\CO_u\= F\cap u^{-1} M(\Z) u.
$$
Let $U$ be the set of all $u\in M(\Q)^\times$ such that
$$
\CO\= F\cap M(\Z)\= F\cap u^{-1} M(\Z)u.
$$
Then $F^\times$ acts on $U$ by multiplication from the
right and $M(\Z)^\times$ acts by multiplication from the
left. It is clear that
$$
\left| M(\Z)^\times \bs U / F^\times\right|\= \left|
M(\Z)^\times \bs \Sigma(\CO)\right|.
$$
So we only have to show that the left hand side equals
$h(\CO)\la_S(\CO)$. For $u\in U$ let
$$
I_u \= F\cap M(\Z)u.
$$
Then $I_u$ is a finitely generated $\CO$-module in $F$. We
claim that the map $\Psi:$
\eqn
M(\Z)^\times \bs U / F^\times &\ra& I(\CO)/F^\times\\
u &\mapsto & I_u,
\endeqn
is surjective and $\la_S(\CO)$ to one. We will show this
through localization and strong approximation. So, for a
prime $p$ let $U_p$ be the set of all $u_p\in
M(\Q_p)^\times$ such that $\CO_p=F_p\cap M(\Z_p)=F_p\cap
u_p^{-1} M(\Z_p)u_p$. We have to show the following:
\begin{itemize}
\item[1.]
For $p\notin S$ the localized map $\Psi_p :
M(\Z_p)^\times\bs U_p / F_p^\times \ra
I(\CO_p)/F_p^\times$ is injective.
\item[2.]
For $p\in S$ the map $\Psi_p$ is $f_p(F)$ to one.
\item[3.]
The map $\Psi$ is surjective.
\end{itemize}
For `1.' let $p\notin S$, let $u_p,v_p\in U_p$, and assume
$$
F_p\cap M(\Z_p)u_p\= F_p\cap M(\Z_p)v_p.
$$
Let $z_p=v_pu_p^{-1}$. Elementary divisor theory implies
that there are $x,y\in M(\Z_p)^\times =
Mat_d(\Z_p)^\times$ such that
$$
z_p\= x{\rm diag}(p^{k_1},\dots,p^{k_d})y,
$$
where diag denotes the diagonal matrix and $k_1\le
k_2\le\cdots\le k_d$ are integers. Replacing $u_p$ by
$yu_p$ and $v_p$ by $x^{-1}v_p$ we may assume that $z$
equals the diagonal matrix. The assumptions then easily
imply $k_1=\cdots =k_d=0$, which gives the first claim.

For `2.' let $p\in S$ and recall that $F_p$ is a local
field and so $h(\CO_p)=1$. So the claim is equivalent to
$$
\left| M(\Z_p)^\times\bs M(\Q_p)^\times /
F_p^\times\right| \= f_p(F).
$$
By Proposition 17.7 of \cite{pierce} it follows that
$$
\left| M(\Z_p)^\times\bs M(\Q_p)^\times\right|
=d=e(M(\Q_p)/\Q_p),
$$
where $e$ denotes the ramification index. If $F_p/Q_p$ is
ramified, then $f_p(F)=1$ and $\left| M(\Z_p)^\times\bs
M(\Q_p)^\times / F_p^\times\right|=1$. If $F_p/\Q_p$ is
unramified , then $f_p(F)=d$ and $\left| M(\Z_p)^\times\bs
M(\Q_p)^\times / F_p^\times\right|=d$ as claimed.

For the surjectivity of $\Psi$ let $I\subset \CO$ be an
ideal. We shall show that there is a $u\in M(\Q)^\times$
such that
$$
F\cap u^{-1} M(\Z)u\= F\cap M(\Z)
$$
and
$$
I\= I_u\= F\cap M(\Z)u.
$$
We shall do this locally. First note that, since $I$ is finitely
generated, there is a finite set of primes $T$ with $T\cap S
=\emptyset$ such that for any $p\notin T\cup S$ the completion
$I_p$ equals $\CO_p$ which is the maximal order of $F_p$. For
these $p$ set $\tilde{u}_p=1$.

Next let $p\in S$. Let $v_p$ be the unique place of $F$ over $p$.
Then $\CO_p=\CO_{v_p}$ is maximal, so it is the valuation ring to
$v_p$ and $I_p=\pi_p^k\CO_p$ for some $k\ge 0$, where $\pi_p$ is
a uniformizing element in $\CO_p$. It follows that at this $p$,
the element $\tilde{u}_p=\pi_p^k Id$ will do the job.

Next let $p\in T$. Then $M(\Z_p)=\Mat_d(\Z_p)$. Let
$\overline{\CO_p} = \CO_p/p\CO_p$ and $\overline{I_p}=I_p/pI_p$.
Then $\overline{\CO_p}$ is a commutative algebra over the field
$\F_p$ with $p$ elements, which implies that
$\overline{\CO_p}\cong\bigoplus_{i=1}^s F_i$, where each $F_i$ is
a finite field extension of $\F_p$. Let $n_i$ be the degree of
$F_i/\F_p$. Then there is an embedding
$\overline{\CO_p}\hookrightarrow \Mat_d(\F_p)$ whose image lies in
$\Mat_{n_1}(\F_p)\times\dots\times\Mat_{n_s}(\F_p)\subset
\Mat_d(\F_p)$. According to the Noether-Skolem Theorem there is a
$\bar{S}\in GL_d(\F_p)$ such that
$\bar{S}\overline{\CO_p}\bar{S}^{-1}\subset
\Mat_{n_1}(\F_p)\times\dots\times\Mat_{n_s}(\F_p)\subset
\Mat_d(\F_p)$. The $\overline{\CO_p}$-ideal $\overline{I_p}$ must
be of the form
$$
\overline{I_p}\= \bigoplus_{i-1}^s \epsilon_i F_i,
$$
where $\epsilon_i\in \{ 0,1\}$. Let $S$ be a matrix in
$GL_d(\Z_p)$ which reduces to $\bar{S}$ modulo $p$ and let
$\tilde{u}_p=S^{-1}(p^{\epsilon_1}Id_{n_1}\times\dots\times
p^{\epsilon_s} Id_{n_s})S$ in $\Mat_d(\Z_p)$. By abuse of
notation we also write $\tilde{u}_p$ for its reduction modulo
$p$. Then we have
$$
\overline{I_p}\= \overline{\CO_p}\cap
\Mat_d(\F_p)\tilde{u}_p.
$$
Let
$$
I_{\tilde{u}_p}=F\cap M(\Z_p)\tilde{u}_p.
$$
Then it follows that
$$
\overline{I_p}\ \cong\ \overline{I_{\tilde{u}_p}}\=
I_{\tilde{u}_p}/pI_{\tilde{u}_p}
$$
and by Theorem 18.6 of \cite{reiner} it follows that $I_p\cong
I_{\tilde{u}_p}$, which implies that there is some $\la\in F_p$
such that $I_p=I_{\tilde{u}_p}\la$. Replacing $\tilde{u}_p$ by
$\tilde{u}_p\la$ and setting $\tilde{u}=(\tilde{u}_p)_p\in
M(\A_{fin})$ we get
$$
I\= F\cap M(\Z)\tilde{u}.
$$
By strong approximation there is a $u\in M(\Q)^\times$ such that
$M(\hat{\Z})u=M(\hat{\Z})\tilde{u}$ and therefore $I=I_u$. \qed

We will summarize the results of this section in the following
proposition.

\begin{proposition}
Let $d$ be a prime $>2$ and let $F/\Q$ be an extension of
degree $d$. Then $F$ embeds into the division algebra
$M(\Q)$ of  degree $d$ if and only if every prime $p$ at
which $M(\Q)$ does not split is non-decomposed in $F$.
Every embedding $\sigma : F\ra M(\Q)$ gives by
intersection with $M(\Z)$ an order $\CO_\sigma$ in $F$.
Every order $\CO$ of $F$, which is maximal at each $p$
where $M(\Q)$ is non-split, occurs in this way. The number
of $M(\Z)^\times$-conjugacy classes of embeddings giving
rise to the same order $\CO$ is equal to
$h(\CO)\la_S(\CO)$.
\end{proposition}

\section{Analysis of the Ruelle zeta function}
From now on we restrict to the case $d=3$. Let $\Ga =\CG(\Z)$;
then $\Ga$ forms a discrete subgroup of $G=\CG(\R)\cong
SL_3(\R)$. Since $M(\Q)$ is a division algebra it follows that
$\CG$ is anisotropic over $\Q$ and so $\Ga$ is cocompact in $G$.
Let $\g_0$ be the real Lie algebra of $G$ and let
$\g=\g_0\otimes_\R\C$ be its complexification. A subgroup
$\Sigma$ of $G$ is called {\it weakly neat} if it is torsion free
and for each $\sigma\in\Sigma$ the adjoint $Ad(\sigma)\in GL(\g)$
does not have a root of unity except $1$ as an eigenvalue. In
other words, a torsion free group $\Sigma$ is weakly neat if for
every $\sigma\in\Sigma$ and every $n\in\N$ the connected
component $G_\sigma^0$ of the centralizer of $\sigma$ coincides
with the connected component $G_{\sigma^n}^0$ of the centralizer
of $\sigma^n$. Further an element $x$ of $G$ is called {\it
regular} if its centralizer is a torus.

\begin{lemma}
The group $\Ga$ is torsion-free and every $\ga\in\Ga$ with
$\ga\ne 1$ is regular. In particular, it follows that $\Ga$ is
weakly neat.
\end{lemma}

\prf Let $\ga\in\Ga$, $\ga\ne 1$; then the centralizer $M(\Q)_\ga$
of $\ga$ in $M(\Q)$ is a division subalgebra of $M(\Q)$ and so it
is either $\Q$, a cubic number field or $M(\Q)$. It cannot be $\Q$
since it contains $\ga$ and if $\ga$ is in $\Q$, it is central so
its centralizer is $M(\Q)$. It can neither be $M(\Q)$, since then
$\ga$ would be in $\Q 1$, say $\ga =r 1$ and then
$1=\det(\ga)=r^3$, thus $\ga=1$ which was excluded. So it follows
that $M(\Q)_\ga$ is a cubic field $F$. Note that this implies that
$\ga$ cannot be a root of unity, since a cubic field does not
contain roots of unity other than $\pm 1$. This shows that $\Ga$
is torsion free. Moreover, we get
$M(\R)_\ga=M(\Q)_\ga\otimes_\Q\R$ so $G_\ga=(F\otimes\R)^1$, the
norm one elements, this is a torus, so $\ga$ is regular. \qed

We have $\g_0 =sl_3(\R)$ and $\g=sl_3(\C)$. Let $b$ be the
Killing form on $\g$. Then
$$
b(X,Y)\= \tr\ad(X)\ad(Y)\= {6}\tr(XY).
$$
Let $K\subset G$ be the maximal compact subgroup $SO(3)$.
Let $\k_0\subset\g_0$ be its Lie algebra and let
$\p_0\subset\g_0$ be the orthogonal space of $\k_0$ with
respect to the form $b$. Then $b$ is positive definite on
$\p_0$ and thus defines a $G$-invariant metric on $X=G/K$,
the symmetric space attached to $G$.

Since $\Ga$ is torsion-free it acts without fixed points on the
contractible space $X$, so $X_\Ga =\Ga\bs X$ is the classifying
space of $\Ga$, in particular, it follows that $\Ga$ is the
fundamental group of $X_\Ga$. We thus obtain a natural bijection
$$
\left\{ \begin{array}{c} {\rm free\ homotopy\ classes}\\
                         {\rm of\ maps}\ S^1\ra X_\Ga\end{array}
                         \right\}
 \leftrightarrow
\left\{ \begin{array}{c} {\rm conjugacy\ classes}\\
                         {[\ga]}\ {\rm in}\ \Ga\end{array}
                         \right\}
$$
For each conjugacy class $[\ga]$ let $X_\ga$ denote the union of
all geodesics lying in the homotopy class $[\ga]$ of $\ga$. It is
known \cite{DKV} that all geodesics in $[\ga]$ have the same
length $l_\ga$ and $X_\ga$ is a submanifold of $X_\Ga$
diffeomorphic to $\Ga_\ga\bs G_\ga /K_\ga$, where $G_\ga$ and
$\Ga_\ga$ are the centralizers of $\ga$ and $K_\ga$ is a maximal
compact subgroup of $G_\ga$.

An element $\ga\in\Ga$ is called {\it primitive} if for
$\sigma\in\Ga$ and $n\in\N$ the equation $\sigma^n=\ga$ implies
that $n=1$. Since every closed geodesic is a positive power of a
unique primitive one it is easy to see that every $\ga\in\Ga$ with
$\ga\ne 1$ is a positive power of a unique primitive element
$\ga_0$. We write $\ga = \ga_0^{\mu(\ga)}$ and call $\mu(\ga)$
the {\it multiplicity} of $\ga$. Clearly primitivity is a property
of the conjugacy class.

Up to conjugacy the group $G$ has two Cartan subgroups,
namely the group of diagonal matrices and the group
$H=AB$, where
$$
A=\left\{
\left. \left( \begin{array}{ccc} a&{}&{}\\ {}&a&{}\\
{}&{}&a^{-2}\end{array}\right)\right| a>0\right\}
$$
and
$$
B=\left(\begin{array}{cc}SO(2)&{}\\ {}&
1\end{array}\right).
$$
Let $P$ denote the parabolic
 $\left( \begin{array}{cc}*&*\\
                         0\ 0&*\end{array}\right)$.
It has a Langlands decomposition $P=MAN$ and $B$ is a compact
Cartan subgroup of
$$
M\cong SL_2^\pm(\R)=\{ x\in Mat_2(\R) | \det(x) =\pm 1\}.
$$
Let
$$
H_1 =\rez{6} \left( \begin{array}{ccc} -1&{}&{}\\ {}& -1&{}\\
{}&{}&2\end{array}\right)\in \a_0 ={\rm Lie}\ A.
$$
Then it follows that $b(H_1)=b(H_1,H_1)=1$. Let $A^-=\{
\exp(tH_1)| t>0\}$ and let $\CE_P(\Ga)$ be the set of all
conjugacy classes $[\ga ]$ in $\Ga$ such that $\ga$ is conjugate
in $G$ to an element $a_\ga b_\ga$ of $A^-B$, and let
$\CE_P^p(\Ga)$ the set of all primitive elements therein. For
$s\in\C$ with $\Re(s)>>0$ let
$$
R_\Ga(s) \= \prod_{[\ga]\in\CE_P^p(\Ga)}(1-e^{-sl_\ga})
$$
be the {\it Ruelle-zeta function} attached to $P$. In
\cite{geom} it is shown that $R_\Ga(s)$ converges for
$\Re(s)>>0$ and that it extends to a meromorphic function
of finite order on the plane. We will show:

\begin{theorem}\label{nust}
The function $R_\Ga(s)$ has a simple zero at $s=1$. Apart
from that, all poles and zeros of $R_\Ga(s)$ are contained
in the union of the interval $[-1,\frac{3}{4}]$ with the
three vertical lines given by $-\rez{2}+i\R$, $i\R$ and
$\rez{2}+i\R$.
\end{theorem}

\prf For any finite dimensional representation $\sigma$ of $M$ let
$$
Z_{P,\sigma}(s) \= \prod_{[\ga]\in\CE_P^p(\Ga)}\prod_{n\ge 0}
\det(1-e^{-sl_\ga} \sigma(b_\ga) S^n(a_\ga b_\ga |\n)),
$$
where $\n =Lie_\C(N)$ and $S^n(a_\ga b_\ga |\n)$ is the $n$-th
symmetric power of the adjoint action of $a_\ga b_\ga$ on $\n$.
In \cite{geom}, Theorem 2.1 it is shown that $Z_{P,\sigma}$
extends to a meromorphic function and that all its poles and
zeros lie in $\R \cup (\rez{2}+i\R)$. Note that $M$ is isomorphic
with $SL_2^\pm(\R)$, the group of real $2\times 2$ matrices of
determinant $\pm 1$. Let $\sigma_0: M\ra GL_2$ denote the standard
representation. In \cite{geom}, Theorem 4.1 it is shown that
$$
R_\Ga(s) \= \frac{Z_{P,1}(s)
Z_{P,1}(s+1)}{Z_{P,\sigma_0}(s+\rez{2})}.
$$
So to complete the proof, it suffices to show the following
proposition.

\begin{proposition}
For $\sigma=1$ the poles and zeros of $Z_{P,\sigma}(s)$ lie in
$[\rez{4},\frac{3}{4}]\cup (\rez{2}+i\R)\cup \{0,1\}$ and the
function $Z_{P,\sigma}(s)$ has a simple zero at $s=1$.

For $\sigma=\sigma_0$ the poles and zeros of $Z_{P,\sigma}$ all
lie in $[0,1]\cup (\rez 2+i\R)$.
\end{proposition}

\prf Let $\hat{G}$ denote the set of all isomorphism
classes of irreducible unitary representations of $G$. The
group $G$ acts on the Hilbert space $L^2(\Ga\bs G)$ by
translations from the right. Since $\Ga\bs G$ is compact
this representation decomposes discretely:
$$
L^2(\Ga\bs G) \= \bigoplus_{\pi\in\hat{G}} N_\Ga(\pi)\pi,
$$
with finite multiplicities $N_\Ga(\pi)$. For $\pi\in\hat{G}$ let
$\pi_K$ denote the $(\g,K)$-module of $K$-finite vectors in
$\pi$. Then the Lie algebra $\n$ acts on $\pi_K$ and we denote by
$H^q(\n,\pi_K)$ the corresponding Lie algebra cohomology
\cite{BorWall}. Let $\m$ denote the complexified Lie algebra of
$M$ and let $\m=\k_M\oplus \p_M$ be its polar decomposition, where
$\k_M$ is the complexified Lie algebra of $K_M=K\cap M$.

In \cite{geom} Theorem 2.1 it is shown that all poles and zeros of
$Z_{P,\sigma}$ lie in $\R\cup (\rez{2} +i\R)$ and that for
$\la\in\a^*$ the (vanishing-) order of $Z_{P,\sigma}$ at
$s=\la(H_1)$ is
$$
\sum_{\pi\in\hat{G}} N_\Ga (\pi) \sum_{p,q\ge 0}
(-1)^{p+q} \dim(H^q(\n,\pi_K)\otimes\wedge^p\p_M\otimes
V_{\breve{\sigma}})_\la^{K_M},
$$
where $(.)_\la$ denotes the generalized $\la$-eigenspace.
Note that the torus $A$ acts trivially on all tensor
factors except $H^q(\n,\pi_K)$. Let $\pi\in\hat{G}$ and
let $\wedge_\pi\in\h^*$ be a representative of its
infinitesimal character. Corollary 3.23 of \cite{HeSch}
says that every character of the $\a$-action on
$H^q(\n,\pi_K)$ is of the form $\mu=w\wedge_\pi +\rho_P$
for some $w\in W(\g,\h)$.

To show the proposition we concentrate on $\sigma$ being trivial
since the other case is similar. Since $\rho(H_1)=-\rez{2}$ we
have to show that for every $\pi$ which has a nonzero
contribution, all eigenvalues $\la$ of $\a$ on $H^p(\n,\pi_k)$
satisfy $-\frac 32\rho_P\le\Re(\la)\le -\rez 2$. By the
isomorphism of $AM$-modules \cite{HeSch}, p. 57:
$$
H_p(\n,\pi_K) \ \cong\ H^{2-p}(\n,\pi_K)\otimes\wedge^2\n
$$
this becomes equivalent to $\rez 2\le\Re(\la)\le \frac 32\rho_P$
whenever $\la$ is an eigenvalue of $\a$ on $H_p(\n,\pi_K)$. Now
fix $\pi\in\hat{G}$ and let $\wedge_\pi\in\h^*$ be a
representative of its infinitesimal character. Then, according to
Corollary 3.32 of \cite{HeSch} we have to show that $-\rez 2
\rho_P\le \Re(\wedge_\pi|_\a)\le\rez 2\rho_P$.

In the case when $\pi$ is induced from the minimal parabolic it
follows that its distributional character $\Theta_\pi$ is zero on
$AB$. By the construction of the test function in \cite{geom}, p.
903, this implies that the contribution of $\pi$ is zero.
Therefore, by the classification of the unitary dual in
\cite{speh}, it remains to consider the case of the trivial
representation and the case when $\pi$ is unitarily induced from
$P=MAN$. So let $\pi=\pi_{\xi,\nu}$ be induced from $P$, where
$\nu$ is imaginary. Then we may assume that $\xi$ is not induced,
since otherwise the double induction formula would lead us back
to the previous case. Let $w\in W$ and $\wedge_\xi$ be the
infinitesimal character of $\xi$. We lift $\wedge_\xi$ to $\d$ by
defining it to be zero on $\a$. Then the Weyl group $W$ will act
on $\wedge_\xi$. We have to show that
$$
-\rez 2 \rho_P\le \Re(w\wedge_\xi|_\a)\le\rez 2\rho_P.
$$
Let us start with $\xi$ being the trivial representation. Then
$\wedge_\xi\matrix t{}{}{-t}$. Lifting $\wedge_\xi$ to $\d$ we get
$\wedge_\xi\left(\begin{array}{ccc}a &&\\ &b& \\ &&
c\end{array}\right)=\rez 2(a-b)$. This vanishes on $A$, so it's
real part is zero. This deals with the case when $w=1$. For $w\in
W$ being the transposition interchanging $b$ and $c$ we get
$w\wedge_\xi\left(\begin{array}{ccc}a &&\\ &b& \\ &&
c\end{array}\right)=\rez 2(a-c)$, which, restricted to $\a$,
coincides with $\rez 2\rho$. This implies the claim for this $w$.
All other Weyl group elements are treated similarly.

It remains to consider the case when $\xi$ is a (limit of)
discrete series representation, so $\xi=\CD_n^+\oplus\CD_n^-$,
where the notation is as in $\cite{knapp}$. Let
$\tau\in\hat{K_M}$ and let $P_\tau :V_\xi\ra V_\xi(\tau)$ be the
projection onto the $\tau$-isotype. For any function $f$ on $G$
which is sufficiently smooth and of sufficient decay the operator
$\pi(f)$ is of trace class. Its trace is
$$
\sum_{\tau\in\hat{K_M}} \int_K\int_{MAN} a^{\nu+\rho}f(k^{-1} man
k)\tr P_\tau \xi(m)P_\tau\ dman\ dk.
$$
Plugging in the test function $f=\Phi$ constructed in \cite{geom},
p.903, this gives
$$
\int_{A^+} a^{\nu+\rho}l_a^{j+1}e^{-sl_a} \tr\xi(f_1) da,
$$
where $f_1$ is the Euler-Poincar\'e function on $M$ attached to
the trivial representation. Then
$$
\tr\xi(f_1)\=\sum_{p=0}^{\dim\p_M} (-1)^p \dim(V_\xi\otimes
\bigwedge^p\p_M)^{K_M}.
$$
Now $K_M\cong O(2)$ and the $K_M$-types can be computed
explicitly. Thus one sees that $\tr\xi(f_1)$ can only be nonzero
if $n=1$ or $n=2$. This means that either $\wedge_\xi\matrix
t{}{}{-t}$ equals $0$ or $t$ respectively, which, in a similar
way to the above, implies the claim. In the case $\sigma=\sigma_0$
the function $f_1$ is replaced by $f_{\sigma_0}$ and one proceeds
in the same fashion. This takes care of all induced
representations. By the classification of the unitary dual of
$GL_3(\R)$ in \cite{speh} it follows that it remains to worry
about the trivial representation only, so let $\pi =triv$ be the
trivial representation of $G$; then the space
$H_0(\n,\pi_K)=\pi_K/\n\pi_K$ is one dimensional with trivial
$\a$-action. This gives a simple zero at $s=1$.
 \qed

\section{Asymptotics of closed geodesics}
For $\ga\in\Ga$ let $N(\ga)=e^{l_\ga}$ and define for
$x>0$:
$$
\pi(x)\= \# \{ [\ga]\in\CE_P^p(\Ga) | N(\ga)\le x\}.
$$
The {\it geodesic prime number theorem} in our context is

\begin{theorem}
For $x\ra\infty$ we have the asymptotic
$$
\pi(x)\ \sim\ \frac{x}{\log x}.
$$
More sharply we have that
$$
\pi(x)\= li(x)+O\left(\frac{x^{\frac{3}{4}}}{\log x}\right)
$$
as $x\ra\infty$ where $li(x)=\int_2^\infty \rez{\log t}dt$ is the
integral logarithm.
\end{theorem}

\prf To simplify the notation in what follows we write $\ga$ for
an element of $\CE_P(\Ga)$ and $\ga_0$ for a primitive element.
If $\ga$ and $\ga_0$ occur in the same formula it is understood
that $\ga_0$ will be the primitive element underlying $\ga$.
Unless otherwise specified, all sums will run either over $\ga$
or $\ga_0$. For $x>0$ let
$$
\psi(x)\= \sum_{N(\ga)\le x} l_{\ga_0}.
$$
For $\Re(s)>1$ we have
\begin{eqnarray*}
\frac{R_\Ga'}{R_\Ga}(s) &=&
\sum_{\ga_0}\frac{l_{\ga_0}e^{-sl_{\ga_0}}}
                                            {1-e^{-sl_{\ga_0}}}\\
&=& \sum_{\ga_0}l_{\ga_0}\sum_{n=1}^\infty e^{-snl_{\ga_0}}\\
&=& \sum_\ga l_{\ga_0}e^{-sl_\ga}.
\end{eqnarray*}
From this point on the argumentation is, up to minor changes the
same as in \cite{randol}, which leads to
$$
\psi(x)\= x+O(x^{\frac{3}{4}}).
$$
From this, the theorem is deduced by standard techniques
\cite{randol} of analytic number theory. \qed

We now finish the proof of the main theorem (\ref{main}). For
this we have to find a division algebra $M(\Q)$ such that for a
given set of primes $S$ with at least two elements we have
$\pi_{S}(x)=\pi(x)$. Firstly, there is a division algebra $M(\Q)$
of degree $3$ such that the set of places at which $M(\Q)$ is
non-split, coincides with the set $S$. This algebra is obtained by
taking the local Brauer-invariants at $p\in S$ to be equal to
$\rez{3}$ or $\frac{2}{3}$ and zero everywhere else in such a way
that they sum to zero in $\Q /\Z$ (\cite{pierce}, Theorem 18.5).
Next, we have a bijection
$$
\CE_P(\Ga)\ra O(S),
$$
given by
$$
[\ga] \mapsto F_\ga\cap M(\Z),
$$
where $F_\ga$ is the centralizer of $\ga$ in $M(\Q)$.
Under this bijection the length $l_\ga$ is transferred to
$\log r(\CO)$. The theorem follows. \qed

\begin{corollary}
We finally note that for $M(\Q)$ chosen as above the
Ruelle zeta function of Theorem \ref{nust} takes the form:
$$
R_\Ga(s/3)\= \prod_{\CO\in O(S)}\left(
1-e^{-sR(\CO)}\right)^{h(\CO)\la_S(\CO)}
$$
\end{corollary}

University of Exeter, Mathematics, Exeter EX4 4QE, Devon, UK


\newpage
\begin{thebibliography}{XXX}


\bibitem{BorWall}
 \bf Borel, A.; Wallach,N.:
 \it Continuous Cohomology, Discrete Groups, and Representations of Reductive Groups.
 \rm Ann. Math. Stud. 94, Princeton 1980.

\bibitem{chand}
 \bf Chandrasekharan, K.:
 \it Introduction to Analytic Number Theory.
 \rm Springer-Verlag 1968.

\bibitem{D-Hitors}
 \bf{Deitmar, A.:}
 \it Higher torsion zeta functions.
 \rm Adv.
Math. 110, 109-128 (1995).

\bibitem{D-Prod}
 \bf{Deitmar, A.:}
 \it Product expansions for zeta functions attached to locally homogeneous spaces.
 \rm Duke Math. J. 82, 71-90 (1996).

\bibitem{D-det}
\bf Deitmar, A.: \it A Determinant Formula for the
generalized Selberg Zeta Function. \rm Quarterly J. Math.
47, 435-453 (1996).

\bibitem{geom}
 \bf Deitmar, A.:
 \it Geometric zeta-functions of locally symmetric spaces.
 \rm Am. J. Math. 122, 887-926 (2000).

\bibitem{DKV}
 \bf Duistermaat, J.J.; Kolk, J.A.C.; Varadarajan, V.S.:
 \it Spectra of locally symmetric manifolds of negative curvature.
 \rm Invent. math. 52 (1979) 27-93.

\bibitem{elstrodt}
 \bf Elstrodt, J.:
 \it Die Selbergsche Spurformel f\"ur kompakte Riemannsche Fl\"achen.
 \rm Jahresbericht der DMV 83, 45-77 (1981).

\bibitem{gauss}
 \bf Gau\ss\ C.F.:
 \it Disquisitiones Arithmeticae.
 \rm Acta Math. 151, 49-151 (1983).

\bibitem{HeSch}
 \bf Hecht, H.; Schmid, W.:
 \it Characters, asymptotics and $\n$-homology of Harish-Chandra modules.
 \rm Acta Math. 151, 49-151 (1983).

\bibitem{hejhal}
 \bf Hejhal, D.:
 \it The Selberg trace formula for $PSL_2(\R)$ I.
 \rm Springer Lecture Notes 548, 1976.

\bibitem{knapp}
 \bf Knapp, A.:
 \it Representation Theory of Semisimple Lie Groups.
 \rm Princeton University Press 1986.

\bibitem{neukirch}
\bf Neukirch, J.: \it Algebraic number theory. \rm
Grundlehren der Mathematischen Wissenschaften. 322.
Berlin: Springer 1999.

\bibitem{pierce}
 \bf Pierce, R. S.:
 \it Associative algebras.
 \rm Springer-Verlag, New York-Berlin, 1982.

\bibitem{randol}
 \bf Randol, B.:
 \it On the asymptotic distribution of closed geodesics on compact Riemann surfaces.
 \rm Trans. AMS 233, 241-247 (1977).

\bibitem{reiner}
 \bf Reiner, I.:
 \it Maximal Orders.
 \rm Academic Press 1975.

\bibitem{sarnak}
 \bf Sarnak, P.:
 \it Class Numbers of Indefinite Binary Quadratic Forms.
 \rm J. Number Theory 15, 229-247 (1982).

\bibitem{siegel}
 \bf Siegel, C. L.:
 \it The average measure of quadratic forms with given determinant and signature.
 \rm Ann. of Math., II. Ser. 45, 667-685 (1944).

\bibitem{speh}
 \bf Speh, B.:
 \it The unitary dual of Gl(3,R) and Gl(4,R).
 \rm Math. Ann. 258, 113-133 (1981).
\end{thebibliography}
\end{document}